\documentclass[11pt]{article}
\usepackage[utf8]{inputenc}
\usepackage{amsmath,amssymb,epsfig,bbm}
\usepackage{stmaryrd}
\usepackage{comment}
\usepackage{color}
\usepackage[T1]{fontenc}

\usepackage[textsize=small]{todonotes}
\usepackage{enumitem}
\usepackage{varwidth}
\setlist{nolistsep}
\usepackage{hyperref}

\newtheorem{defi}{Definition}[section]
\newtheorem{prop}[defi]{Proposition}
\newtheorem{theo}[defi]{Theorem}
\newtheorem{conj}[defi]{Conjecture}
\newtheorem{lemm}[defi]{Lemma}
\newtheorem{coro}[defi]{Corollary}
\newtheorem{rema}[defi]{Remark}
\newtheorem{exem}[defi]{Example}
\newtheorem{exems}[defi]{Examples}

\newcommand{\bdefi}{\begin{defi}}
\newcommand{\edefi}{\end{defi}}
\newcommand{\bprop}{\begin{prop}}
\newcommand{\eprop}{\end{prop}}
\newcommand{\btheo}{\begin{theo}}
\newcommand{\etheo}{\end{theo}}
\newcommand{\blemm}{\begin{lemm}}
\newcommand{\brema}{\begin{rema}}
\newcommand{\erema}{\end{rema}}
\newcommand{\bexer}{\begin{exem}}
\newcommand{\eexer}{\end{exem}}
\newcommand{\bexems}{\begin{exems}}
\newcommand{\eexems}{\end{exems}}
\newcommand{\bconj}{\begin{conj}}
\newcommand{\econj}{\end{conj}}
\newcommand{\elemm}{\end{lemm}}
\newcommand{\bcoro}{\begin{coro}}
\newcommand{\ecoro}{\end{coro}}
\newcommand{\dem}{\noindent{\bf Proof. }}

\usepackage{mathrsfs}
\renewcommand\mathcal{\mathscr}

\newcommand{\F}{{\cal F}}

\newcommand{\M}{{\cal M}}

\newcommand{\OOO}{{\cal O}}

\newcommand{\maths}[1]{{\mathbb #1}}  

\newcommand{\CC}{\maths{C}}

\newcommand{\HH}{\maths{H}}

\newcommand{\NN}{\maths{N}}

\newcommand{\PP}{\maths{P}}
\newcommand{\QQ}{\maths{Q}}
\newcommand{\RR}{\maths{R}}
\newcommand{\SSS}{\maths{S}}

\newcommand{\UU}{\maths{U}}

\newcommand{\ZZ}{\maths{Z}}

\newcommand{\ra}{\rightarrow}
\newcommand{\bs}{\backslash}

\newcommand{\wt}[1]{{\widetilde{#1}}}

\newcommand{\ga}{\gamma}
\newcommand{\Ga}{\Gamma}

\newcommand{\cqfd}{\hfill$\Box$}

\renewcommand{\Re}{{\operatorname{Re}}}
\renewcommand{\Im}{{\operatorname{Im}}}

\newcommand{\PSL}{\operatorname{PSL}}
\newcommand{\SL}{\operatorname{SL}}
\newcommand{\GL}{\operatorname{GL}}

\newcommand{\Heis}{\operatorname{Heis}}

\newcommand{\HS}{\mathcal{H\!S}}
\newcommand{\mat}{\M}
\newcommand{\Gal}{\operatorname{Gal}}

\newcommand{\stab}{\operatorname{Stab}}

\newcommand{\maxcf}{\operatorname{\F_\CC}}

\newcommand{\hdc}{{\HH}^2_\CC}

\newcommand{\PSLOK}{\operatorname{PSL}_{2}(\OOO_K)}

\newcommand{\PSLC}{\operatorname{PSL}_{2}(\CC)}
\newcommand{\SLR}{\operatorname{SL}_{2}(\RR)}
\newcommand{\PSLR}{\operatorname{PSL}_{2}(\RR)}

\newcommand{\SU}{\operatorname{SU}}
\newcommand{\PSU}{\operatorname{PSU}}

\newcommand{\tr}{\operatorname{\tt tr}}

\newcounter{fig}

\title{A classification of $\CC$-Fuchsian subgroups \\
of Picard modular groups}
\author{Jouni Parkkonen \and Fr\'ed\'eric Paulin} 

\begin{document}
\bibliographystyle{../alphanum}
\maketitle
\begin{abstract}
 Given an imaginary quadratic extension $K$ of $\QQ$, we give a
  classification of the maximal nonelementary subgroups of the Picard
  modular group $\PSU_{1,2}(\OOO_K)$ preserving a complex geodesic in
  the complex hyperbolic plane $\HH^2_\CC$. Complementing work of
  Holzapfel, Chinburg-Stover and M\"oller-Toledo, we show that these
  maximal $\CC$-Fuchsian subgroups are arithmetic, arising from a
  quaternion algebra $\Big(\!\begin{array}{c}
    D\,,D_K\\\hline\QQ\end{array} \!\Big)$ for some explicit $D\in\NN-\{0\}$
  and $D_K$ the discriminant of $K$.  We thus prove the existence of
  infinitely many orbits of $K$-arithmetic chains in the hypersphere
  of $\PP_2(\CC)$. 
  
\end{abstract}

\section{Introduction}

Let $h$ be a Hermitian form with signature $(1,2)$ on $\CC^3$. The
projective special unitary subgroup $\PSU_h$ of $h$ contains two
conjugacy classes of Lie subgroups isomorphic to $\PSLR$. The
subgroups in one class preserve a complex projective line for the
projective action of $\PSU_h$ on the projective plane $\PP_2(\CC)$,
and those of the other class preserve a totally real subspace.  The
groups $\PSLR$ and $\PSU_h$ act as the groups of orientation
preserving isometries, respectively, on the upper halfplane model of
the real hyperbolic space and on the projective model of the complex
hyperbolic plane defined using the form $h$.  If $\Ga$ is a discrete
subgroup of $\PSU_h$, the intersections of $\Ga$ with the Lie
subgroups isomorphic to $\PSLR$ are its {\em Fuchsian subgroups} and
the Fuchsian subgroups preserving a complex projective line are called
{\it $\CC$-Fuchsian subgroups}.  We refer to Section \ref{sec:cxhyp}
for more precise definitions and comments on the terminology.

Let $K$ be an imaginary quadratic number field, with discriminant
$D_K$ and ring of integers $\OOO_K$. We consider the Hermitian form
$h$ defined by
$$
(z_0,z_1,z_2)\mapsto 
-z_0\,\overline{z_2}-z_2\,\overline{z_0}+ z_1\overline{z_1}\;.
$$ 
The {\it Picard modular group} $\Ga_K=\PSU_h(\OOO_K)$ is a nonuniform
arithmetic lattice of $\PSU_h$, see for instance
\cite[Chap.~5]{Holzapfel98} and subsequent works of Falbel, Parker,
Francsics, Lax, Xie-Wang-Jiang, and many others, for information on
these groups.  In this paper, we classify the maximal $\CC$-Fuchsian
subgroups of $\Ga_K$, and we explicit their arithmetic structures.

When $G=\PSLC$, there is exactly one conjugacy class of Lie subgroups
of $G$ isomorphic to $\PSLR$.  When $\Ga$ is the {\it Bianchi group}
$\PSLOK$, the analogous classification is due to Maclachlan and Reid
(see \cite{Maclachlan86, MacRei91} and \cite[Chap.~9]{MacRei03}). They
proved that the maximal nonelementary Fuchsian subgroups of $\PSLOK$
are commensurable up to conjugacy with the stabilisers of the circles
$|z|^2=D$ for $D\in\NN-\{0\}$, when $\PSLC$ acts projectively (by
homographies) on the projective line $\PP_1(\CC)= \CC \cup\{\infty\}$,
and that all these subgroups arise from explicit quaternion algebras
over $\QQ$.  For information on Bianchi groups, see for instance
\cite{Fine89} and the references of \cite{MacRei91}.

More generally, given a semisimple connected real Lie group $G$ with
finite center and without compact factor, there is a nonempty finite
set of infinite conjugacy classes of Lie subgroups of $G$ locally
isomorphic to $\SLR$, unless $G$ itself is locally isomorphic to
$\SLR$. The structure of the set of these subgroups plays an important
role for the classification of the linear representations of $G$, and
for the classification of the groups $G$ themselves, see for instance
\cite{Knapp05, Serre92} among others. Given a discrete subgroup $\Ga$
of $G$, it is again interesting to study the {\it Fuchsian subgroups}
of $\Ga$, that is, the intersections of $\Ga$ with these Lie
subgroups, to classify the maximal ones and to see, when $\Ga$ is
arithmetic, if its maximal Fuchsian subgroups are also arithmetic (see
Proposition \ref{prop:arithmeticgeneral} for a positive answer) with
an explicit arithmetic structure. From now on, $G=\PSU_h$.

We first prove (see Proposition \ref{prop:Fuchsianrationalpoint} and
just after) that a nonelementary $\CC$-Fuchsian subgroup $\Ga'$ of
$\Ga_K$ preserves a unique projective point $[z_0:z_1:z_2]$ with
$z_0,z_1,z_2$ relatively prime in $\OOO_K$. 
We define the {\it discriminant} of $\Ga'$ as
$\Delta_{\Ga'}=h(z_0,z_1,z_2)$.
%
%
For any positive natural number $D$, let
$$
\Ga_{K,D}=\stab_{\Ga_K}[-D:0:1]\,.
$$
In Section \ref{sec:class}, we prove the following classification
result (see \cite[Thm.~1]{MacRei91} and \cite[Thm.~9.6.2]{MacRei03} in
the Bianchi group case).

\btheo\label{theo:classification} Let $D\in\NN-\{0\}$. The set of
$\Ga_K$-conjugacy classes of maximal nonelementary $\CC$-Fuchsian
subgroups of $\Ga_K$ with discriminant $D$ is finite and nonzero. Every
maximal nonelementary $\CC$-Fuchsian subgroup of $\Ga_K$ with
discriminant $D$ is commensurable up to conjugacy in $\PSU_h$ with
$\Ga_{K,2D}$.  
\etheo

In the course of the proof of this result, we prove a criterion for
when two groups $\Ga_{K,D}$ for $D\in\NN-\{0\}$ are commensurable up
to conjugacy in $\PSU_h$. A further application of this condition
shows that every maximal nonelementary $\CC$-Fuchsian subgroup of
$\Ga_K$ is commensurable up to conjugacy in $\PSU_h$ with $\Ga_{K,D}$
for a squarefree natural number $D$.

\smallskip
%
%
Recall (see for instance \cite{Goldman99}) that a {\it chain}
\footnote{a notion attributed to von Staudt in \cite[footnote
  3]{Cartan32}} is the intersection of the {\it Poincaré
  hypersphere} $$\HS=\{[z]\in\PP_2(\CC)\;:\;h(z)=0\}$$ with a complex
projective line (if nonempty and not a singleton). It is {\it
  $K$-arithmetic} if its stabiliser in $\Ga_K$ has a dense orbit in
it.  

\bcoro \label{coro:intro} There are infinitely many $\Ga_K$-orbits of
$K$-arithmetic chains in the hypersphere $\HS$.  
\ecoro

The figure below shows part of the image under vertical projection in
the Heisenberg group of the orbit under $\Ga_K$ of a $K$-arithmetic
chain whose stabiliser has discriminant $10$, when $K=\QQ[i]$.

\begin{center}
\includegraphics[width=10cm]{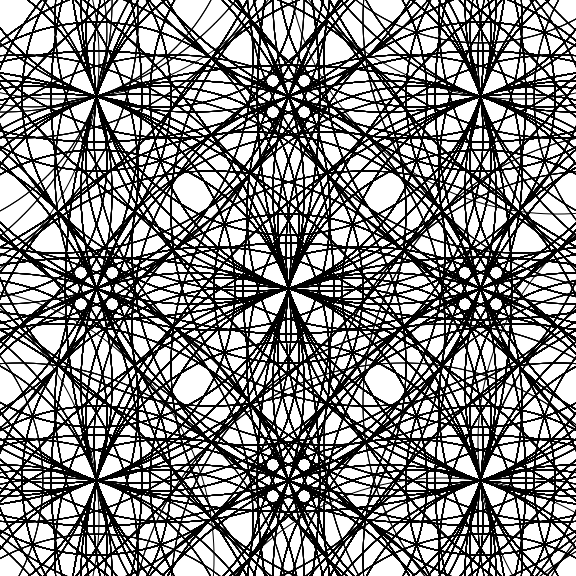}
\end{center}

We say that a subgroup of $\PSU_h$ {\em arises from a quaternion
  algebra} $A$ defined over $\QQ$ if it is commensurable with
$\sigma(A(\ZZ)^1)$ for some $\QQ$-algebra morphism $\sigma:
A\to\mat_3(\CC)$.  In Section \ref{sec:quaternion}, we prove the
following result (see \cite[Thm.~9.6.3]{MacRei03} in the Bianchi group
case).

\btheo\label{theo:quaternion} Every nonelementary $\CC$-Fuchsian
subgroup of $\Ga_K$ of discriminant $D$ is conjugate in $\PSU_h$ to a
subgroup of $\PSU_h$ arising from the quaternion algebra
$\big(\frac{D,\,D_K}{\QQ} \big)$.  
\etheo

The classification of the quaternion algebras over $\QQ$ then allow to
classify up to commensurability and conjugacy the maximal
nonelementary $\CC$-Fuchsian subgroups of $\Ga_K$: two such groups,
with discriminant $D$ and $D'$ are commensurable up to conjugacy if
and only if the quaternion algebras $\big(\frac{D,\,D_K}{\QQ} \big)$
and $\big(\frac{D',\,D_K}{\QQ} \big)$ are isomorphic. This holds for
instance if and only if the quadratic forms $D_Kx^2+Dy^2-DD_Kz^2$ and
$D_Kx^2+D'y^2-D'D_Kz^2$ are equivalent over $\QQ$.

As was mentioned to us by M.~Stover after we posted a first version
of this paper on ArXiv, the existence of a bijection between wide
commensurability classes of $\CC$-Fuchsian subgroups of $\Ga_K$ and
isomorphism classes of quaternion algebras over $\QQ$ unramified at
infinity and ramified at all finite places which do not split in
$K/\QQ$ is a particular case of the 2011 unpublished preprint
\cite{ChiSto11} (see Theorem 2.2 in its version 3), which proves such
a result for all arithmetic lattices of simple type in $\SU_{2,1}$. In
particular, the existence of this bijection (and our Corollary
\ref{coro:manyclasses}) should be attributed to
Chinburg-Stover (although
they say it was known by experts). Möller-Toledo in \cite{MolTol14} (a reference we were
also not aware of for the first draft of this paper) also give a
description of the quotients by the maximal $\CC$-Fuchsian subgroups
of the real hyperbolic planes they preserve, and more generally of all
Shimura curves in Shimura surfaces of the first type. We believe
that our precise correspondence brings interesting effective and
geometric informations.

\medskip\noindent{\small {\it Acknowledgements: } The first author
  thanks the V\"ais\"al\"a foundation and the FIM of ETH Z\"urich for
  their support during the preparation of this paper. The second
  author thanks the Väisälä foundation and its financial support for a
  fruitful visit to the University of Jyväskylä and the nordic
  snows. This work is supported by the NSF Grant no 093207800, while
  the second author was in residence at the MSRI, Berkeley CA, during
  the Spring 2015 semester. We thank Y.~Benoist and M.~Burger for
  interesting discussions on this paper. We warmly thank M.~Stover for
  informing us about the paper \cite{ChiSto11} and many other
  references, including \cite{MolTol14}.}

\section{The complex hyperbolic plane}
\label{sec:cxhyp}

\medskip
Let $h$ be the nondegenerate Hermitian form 
$$
h(z)=-z_0\,\overline{z_2}-z_2\,\overline{z_0} + |z_1|^2=
-2\,\Re(z_0\,\overline{z_2})+ |z_1|^2
$$ 
of signature $(1,2)$ on $\CC^3$ with coordinates $(z_0,z_1,z_2)$, and
let $\langle\cdot,\cdot\rangle$ be the associated Hermitian product.
The point $z=(z_0,z_1,z_2)\in\CC^3$ and the corresponding element
$[z]=[z_0:z_1:z_2]\in \PP_2(\CC)$ (using homogeneous coordinates) is
{\em negative, null or positive} according to whether $h(z)<0$,
$h(z)=0$ or $h(z)>0$.  The {\it negative/null/positive cone} of $h$ is
the subset of negative/null/positive elements of $ \PP_2(\CC)$.

The negative cone of $h$ endowed with the distance $d$ defined by
$$
\cosh^2d([z],[w])=\frac{|\langle z,w\rangle|^2}{h(z)\,h(w)}
$$
is the {\em complex hyperbolic plane} $\hdc$. The distance $d$ is the
distance of a Riemannian metric with pinched negative sectional
curvature $-4\le K\le -1$.  The linear action of the special unitary
group of $h$
$$
\SU_h=\{g\in \SL_{3}(\CC)\;:\; h\circ g=h\}
$$
on $\CC^{3}$ induces a projective action on $\PP_2(\CC)$ with kernel
$\UU_{3}\operatorname{Id}$, where $\UU_{3}$ is the group of third
roots of unity.  This action preserves the negative, null and positive
cones of $h$, and is transitive on each of them. The restriction to
$\hdc$ of the quotient group $\PSU_h= \SU_h/(\UU_{3}
\operatorname{Id})$ of $\SU_h$ is the orientation-preserving isometry
group of $\hdc$.

The null cone of $h$ is the {\em Poincar\'e hypersphere} $\HS$, which
is naturally identified with the boundary at infinity of $\hdc$.  The
{\em Heisenberg group}
$$
\Heis_3=\{[w_0:w:1]\in\CC\times\CC:2\,\Re\, w_0=|w|^2\} 
$$
acts isometrically on $\hdc$ and simply transitively on $\HS-
\{[1:0:0]\}$ by the action induced by the matrix representation
$$
[w_0:w:1]\mapsto \begin{pmatrix} 1&\overline{w}& w_0\\0&1&w \\0&0&1
\end{pmatrix}
$$
of $\Heis_3$ in $\SU_h$. The projective transformations induced by
these matrices are called {\em Heisenberg translations}.

If a complex projective line meets $\hdc$, its intersection with
$\hdc$ is a totally geodesic submanifold of $\hdc$, called a {\em
  complex geodesic}. The intersection of a complex projective line in
$\PP_2(\CC)$ with the Poincaré hypersphere is called a {\em chain}, if
nonempty and not reduced to a point.  Each complex projective line $L$
in $\PP_2(\CC)$ meeting $\hdc$ (or its associated complex geodesic
$L\cap\hdc$, or its associated chain $L\cap\HS$) is {\em polar} to a
unique positive point $P_L\in\PP_2(\CC)$, that is, $\langle z,
P_L\rangle=0$ for all $z\in L$ (or equivalently $z\in L\cap \hdc$ or
$z\in L\cap \HS$).  This element $P_L$ is the {\em polar point} of the
projective line $L$, of the complex geodesic $L\cap\hdc$ and of the
chain $L\cap\HS$. Conversely, for each positive point $P$, there is a
unique complex projective line $P^\perp$ polar to $P$, the {\em polar
  line} of $P$. The intersection of $P^\perp$ with $\hdc$ is a complex
geodesic.

An easy computation (using for instance Equation (42) in
\cite{ParPau10GT}) shows that
\begin{equation}\label{eq:calcstabpospoint}
  \stab_{\SU_h}[0:1:0]=\Big\{
\begin{pmatrix} \zeta a&0&i\zeta b\\0&\zeta^{-2}&0\\
-i\zeta c&0&\zeta d\end{pmatrix}
\;:\; \begin{array}{l}a,b,c,d\in\RR,\;\zeta\in\CC\\
ad-bc=1,\;|\zeta| =1\end{array}\Big\}\,.
\end{equation}
In particular, $\stab_{\SU_h}[0:1:0]$ is isomorphic to $\SSS^1\times
\PSLR$, and $\stab_{\PSU_h}[0:1:0]$ is also isomorphic to $\SSS^1
\times \PSLR$.  More generally, if $P=[z_0:z_1:z_2]$ is a positive
point in $\PP_2(\CC)$, then by Equation \eqref{eq:calcstabpospoint},
its stabiliser in $\PSU_h$ is the direct product of a Lie group
embedding of $\PSLR$ in $\PSU_h$ preserving the complex geodesic polar
to $P$, with the group of complex reflections with fixed point set the
projective line polar to $P$.

The polar chain of $P$ is
$$
C_P=
\{[w_0:w_1:w_2]\in\PP_2(\CC)\;:\; h(w_0,w_1,w_2)=
\langle(w_0,w_1,w_2),(z_0,z_1,z_2)\rangle=0\}\,,
$$ 
that is $C_P\cap \Heis_3$ is the set of $[w_0:w:1]\in\Heis_3$
satisfying the equation
$$
\big(\frac{|w|^2}{2}+i\,\Im\;w_0\big)\,\overline{z_2}-w\,\overline{z_1}
+\overline{z_0}=0\;.
$$
When $z_2\neq 0$, in the coordinates $(w,2\,\Im\;w_0) \in
\CC\times\RR$ of $[w_0:w:1]\in\Heis_3$, this is the equation of an
ellipse, whose image under the {\em vertical projection} $[w_0:w:1]
\mapsto w$ is the circle with center $\frac{\overline{z_1}}
{\overline{z_2}}$ and radius $\frac{\sqrt{h(z_0,z_1,z_2)}}{|z_2|}$ in
$\CC$ given by the equation
$$
|w|^2-2\,\Re\big(w\,\frac{\overline{z_1}}{\overline{z_2}}\big)
+2\,\Re\big(\, \frac{\overline{z_0}}{\overline{z_2}}\,\big)=0\,.
$$
If $z_2=0$, then
$C_P\cap \Heis_3$ is the vertical affine line over
$\frac{\overline{z_1}}{\overline{z_2}}$.

We refer to Goldman \cite[p.~67]{Goldman99} and Parker \cite{Parker10}
for the basic properties of $\hdc$. These references use different
Hermitian forms of signature $(1,2)$ to define the complex hyperbolic
plane, and the curvature is often normalised differently from our
definitions.  Our choices are consistent with \cite{ParPau10GT} and
\cite{ParPauHeis}.

\section{Classification of $\CC$-Fuchsian subgroups 
of $\Ga_K$}
\label{sec:class}

Before starting to study Fuchsian subgroups of discrete subgroups of
$\PSU_h$, let us mention that it is a very general fact that the
maximal {\it nonelementary} (that is, not virtually cyclic) Fuchsian
subgroups of arithmetic subgroups of $\PSU_h$ are automatically
(arithmetic) lattices of the copy of $\PSLR$ containing them.

\bprop \label{prop:arithmeticgeneral} 
Let $G$ be a semisimple connected real Lie group with finite center
and without compact factor, and let $\Ga$ be a maximal nonelementary
Fuchsian subgroup of an arithmetic subgroup $\wt\Ga$ of $G$. Then
$\Ga$ is an arithmetic lattice in the copy of the group locally
isomorphic to $\SLR$ containing it.  
\eprop

One of the main points of what follows will be to determine
explicitly the arithmetic structure of $\Ga$, that is the
$\QQ$-structure thus constructed on the group locally isomorphic to
$\SLR$ containing it, relating it to the arithmetic structure of
$\wt\Ga$, that is the given $\QQ$-structure on $G$.

\medskip \dem Let $\underline{G}$ be a semisimple connected algebraic
group defined over $\QQ$, let $\underline{H}$ be an algebraic subgroup
of $G$ defined over $\RR$ locally isomorphic to $\SL_2$, and assume
that $\Ga=\underline{H}(\RR)\cap G(\ZZ)$ is nonelementary in
$\underline{H}(\RR)$.  As a nonelementary subgroup of a group locally
isomorphic to $\SL_2$ is Zariski-dense in it, and as the Zariski
closure of a subgroup of $G(\ZZ)$ is defined over $\QQ$, we hence have
that $\underline{H}$ is defined over $\QQ$. Therefore by the
Borel-Harish-Chandra theorem \cite[Thm.~7.8]{BorelHarishChandra62},
$\Ga= \underline{H}(\ZZ)$ is an arithmetic lattice in
$\underline{H}(\RR)$.  Since the copies of subgroups of $G$ locally
isomorphic to $\SLR$ are algebraic, the result follows.  
\cqfd

\bigskip
Let $K$ be an imaginary quadratic number field, with $D_K$ its
discriminant, $\OOO_K$  its ring of integers, $\tr:z\mapsto
z+\overline{z}$ its trace and $N:z\mapsto|z|^2$ its norm.  The {\em
  Picard modular group} of $K$, that we denote by $\Ga_K=
\PSU_h(\OOO_K)$, consists of the images in $\PSU_h$ of matrices of
$\SU_h$ with coefficients in $\OOO_K$.  It is a nonuniform arithmetic
lattice by the result of Borel and Harish-Chandra cited above.

A discrete subgroup $\Ga$ of $\PSU_h$ is an {\em extended
  $\CC$-Fuchsian subgroup} if it satisfies one of the following
equivalent conditions
\begin{enumerate}
\item $\Ga$ preserves a complex projective line of $\PP_2(\CC)$
  meeting $\hdc$,
\item $\Ga$ fixes a positive point in $\PP_2(\CC)$,
\item $\Ga$ preserves a chain.
\end{enumerate}
Many references, see for example \cite{FalPar03}, do not use the word
``extended''. But as defined in the introduction, in this paper, a
{\it $\CC$-Fuchsian subgroup} is a discrete subgroup of $\PSU_h$
preserving a complex geodesic in $\hdc$ and inducing the parallel
transport on its unit normal bundle. It is the image of a {\it
  Fuchsian group} (that is, a discrete subgroup of $\PSLR$) by a Lie
group embedding of $\PSLR$ in $\PSU_h$. The extended $\CC$-Fuchsian
subgroups are then finite extensions of $\CC$-Fuchsian subgroups by
finite groups of complex reflections fixing the projective line or
positive point or chain in the definition above. In particular, up to
commensurability, the notions of extended $\CC$-Fuchsian subgroups and
of $\CC$-Fuchsian subgroups coincide. The $\CC$-Fuchsian lattices have
been studied under a different viewpoint than our differential
geometric one, as fundamental groups of arithmetic curves on ball
quotient surfaces or Shimura curves in Shimura surfaces, by many
authors, see for instance \cite{Kudla78, Holzapfel83, Holzapfel98,
  MolTol14} and their references.

\medskip An element of $\Ga_K$ is {\it $K$-irreducible} if it does not
preserve a point or a line defined over $K$ in $\PP_2(\CC)$.  An
element of $\PP_2(\CC)$ is {\em rational} if it lies in $\PP_2(K)$. Note
that the polar line of a positive rational point of $\PP_2(\CC)$ is
defined over $K$. The group $\PSU_h(K)$, image of $\SU_h(K)= \SU_h\cap
\SL_3(K)$ in $\PSU_h$, preserves $\PP_2(K)$, but in general acts
transitively on neither the positive, the null nor the negative points
of $\PP_2(K)$.

The Galois group $\Gal(\CC|K)$ acts on $\PP_2(\CC)$ by $\sigma
[z_0:z_1:z_2]= [\sigma z_0:\sigma z_1:\sigma z_2]$, and fixes
 $\PP_2(K)$ pointwise. A positive point $z\in \PP_2(\CC)$ is {\it
  Hermitian cubic} over $K$ if it is cubic over $K$ (that is, if its
orbit under $\Gal(\CC|K)$ has exactly three points), and if its other
Galois conjugates $z',z''$ over $K$ are null elements in the polar
line of $z$.

The following result, analog to \cite[Prop.~9.6.1]{MacRei03} in the
Bianchi group case, strengthens one direction of
\cite[Lem.~1.2]{MolTol14}.

\bprop\label{prop:Fuchsianrationalpoint} A nonelementary extended
$\CC$-Fuchsian subgroup $\Ga$ of $\Ga_K$ fixes a unique rational point
in $\PP_2(\CC)$. This point is positive and it is the polar point of
the unique complex geodesic preserved by $\Ga$. 
\eprop

\dem If $\alpha\in\PSU_h$ is loxodromic, let $\alpha_-,\alpha_+
\in \partial_\infty\hdc$ be its repelling and attracting fixed points,
and let $\alpha_0$ be its positive fixed point. Since the two
projective lines tangent to the hypersphere $\HS$ at $\alpha_-$ and
$\alpha_+$ are invariant under $\alpha$, their unique intersection
point is fixed by $\Ga$, therefore is equal to $\alpha_0$. In particular,
$\alpha_0$ is polar to the complex projective line through
$\alpha_-,\alpha_+$ (see also \cite[Lemma 6.6]{Parker10} for a more
analytic proof).

Let $L$ be the complex projective line preserved by $\Ga$, which meets
$\hdc$.  As $\Ga$ is not elementary, there are loxodromic elements
$\alpha, \beta \in\Ga$ such that their sets of fixed points in
$\partial_\infty\hdc \cap L$ are disjoint. Since $L$ passes through
$\alpha_-,\alpha_+$ as well as through $\beta_-,\beta_+$, and by the
uniqueness of the polar point to $L$, we hence have $\alpha_0=
\beta_0$.

As $\alpha$ and $\beta$ have infinite order, one of them cannot be
$K$-irreducible. Otherwise, if both were $K$-irreducible, then by
\cite[Prop.~18]{ParPauHeis}, the point $\alpha_0=\beta_0$ would be
Hermitian cubic and its orbit under $\Gal(\CC|K)$ would be
$\{\alpha_-,\alpha_+,\alpha_0\} = \{\beta_-,\beta_+,\beta_0\}$, a
contradiction.  Assume then for instance that $\alpha$ preserves a
line or a point defined over $K$.  As any projective subspace
preserved by $\alpha$ is a combination of $\alpha_-$, $\alpha_+$ and
$\alpha_0$, and as $\alpha_-$ and $\alpha_+$ are not defined over $K$,
it follows that $\alpha_0$ is rational.  
\cqfd

\medskip 
Let $\Ga$ be a nonelementary extended $\CC$-Fuchsian subgroup of
$\Ga_K$. By the previous proposition, $\Ga$ fixes a unique rational
point $P_\Ga$ in $\PP_2(\CC)$, which may be written $P_\Ga=
[z_0:z_1:z_2]$ with $z_0,z_1,z_2\in\OOO_K$ relatively prime. Such a
writing is unique up to the simultaneous multiplication of
$z_0,z_1,z_2$ by a unit in $\OOO_K$.  Since the units in $\OOO_K$ have
norm $1$, and since the trace and norm of $K$ take integral values on
the integers of $K$, the number
$$
\Delta_\Ga=h(z_0,z_1,z_2)=N(z_1)-\tr(z_0\,\overline{z_2})\in\ZZ
$$
is well defined, we call it the {\it discriminant} of $\Ga$.  As $P$
is positive, we have $\Delta_\Ga\in\NN-\{0\}$. The radius of the
vertical projection of the polar chain of $P_\Ga$ is hence
$\frac{\sqrt{\Delta_\Ga}}{|z_2|}$. The discriminant of $\Ga$ depends
only on the conjugacy class of $\Ga$ in $\Ga_K$: for every
$\ga\in\Ga_K$, since by uniqueness we have $P_{\ga\Ga\ga^{-1}}= \ga
P_\Ga$, we have
$$
\Delta_{\ga\Ga\ga^{-1}}=\Delta_\Ga\;.
$$

\medskip 
A chain $C$ is {\em ($K$-)arithmetic} if its stabiliser in $\Ga_K$ has
a dense orbit in $C$. The following result along with Proposition
\ref{prop:Fuchsianrationalpoint} justifies this terminology. This result is well known, and it is the other direction of
\cite[Lem.~1.2]{MolTol14}, see also \cite[Prop.~1.5,\S
III.1]{Holzapfel83} and \cite[\S 3]{Kudla78}.  We give a proof, which
is a bit different, for the sake of completeness.

\bprop \label{prop:stabmaxextCfuchs} 
The stabiliser $\stab_{\Ga_K}P$ of any positive rational point
$P\in \PP_2(K)$ is a maximal nonelementary extended $\CC$-Fuchsian
subgroup of $\Ga_K$, whose invariant chain is arithmetic.  
\eprop

\dem  Let $\underline{G}$
be the linear algebraic group defined over $\QQ$, such that
$\underline{G}(\ZZ)=\PSU_h(\OOO_K)$ and $\underline{G}(\RR) = \PSU_h$.
We endow $\PP_2(\CC)$ with the $\QQ$-structure $\underline X$ whose
$\QQ$-points are $\PP_2(K)$ so that the action of $\underline G$ on
$\underline X$ is defined over $\QQ$.

As seen in Section \ref{sec:cxhyp}, the set of real points of
$\stab_{\underline G}P$ is isomorphic to $\SSS^1\times \PSLR$ as a
real Lie group. The group $\stab_{\underline G}P$ is reductive and it
has a (semisimple) Levi subgroup $\underline H$ defined over $\QQ$,
such that $\underline H(\RR)$ is isomorphic to $\PSLR$. By a theorem
of Borel-Harish-Chandra \cite{BorelHarishChandra62}, the group
$\underline H(\ZZ)$ is an arithmetic lattice in $\underline H(\RR)$,
which (preserves the projective line polar to $P$ and) is contained in
$\stab_{\Ga_K}P$.  As $\underline H(\ZZ)$ is a lattice in $\underline
H(\RR)$, the group $\stab_{\Ga_K}P$ is nonelementary and has a dense
orbit in the chain $P^\perp\cap\HS$.  
\cqfd

\medskip Recall that in the coordinates $(w,-2\,\Im\; w_0)$ of
$\Heis_3$, the chains are ellipses whose images under the vertical
projection are Euclidean circles (see also \cite[\S 4.3]{Goldman99}).
The figure in the introduction is the vertical projection of part of
the orbit under $\Ga_K$ of the chain $[-5:0:1]^\perp\cap\HS$ when
$K=\QQ[i]$, so that $\Ga_K$ is the Gauss-Picard modular group, whose
generators have been given by \cite{FalFraPar11}. The figure shows the
square $|\,\Re\, z|, |\,\Im \,z|\le 1.5$ in $\CC$ with projections of
chains whose diameter is at least $1$.  In the figures below,
$K=\QQ[\omega]$, where $\omega$ is a primitive third root of unity, so
that $\Ga_K$ is the Eisenstein-Picard modular group, whose generators
have been given by \cite{FalPar06}.  The first figure shows part of the
orbit of $[-1:0:1]^\perp\cap\HS$ and the second figure shows part of
the orbit of $[-2:0:1]^\perp\cap\HS$.  They both show the square
$|\,\Re\, z|, |\,\Im\, z|\le 1$ in $\CC$ with projections of chains
whose diameter is at least $0.5$ in the first figure and at least 0.75 in the second.

\begin{center}
\includegraphics[width=10cm]{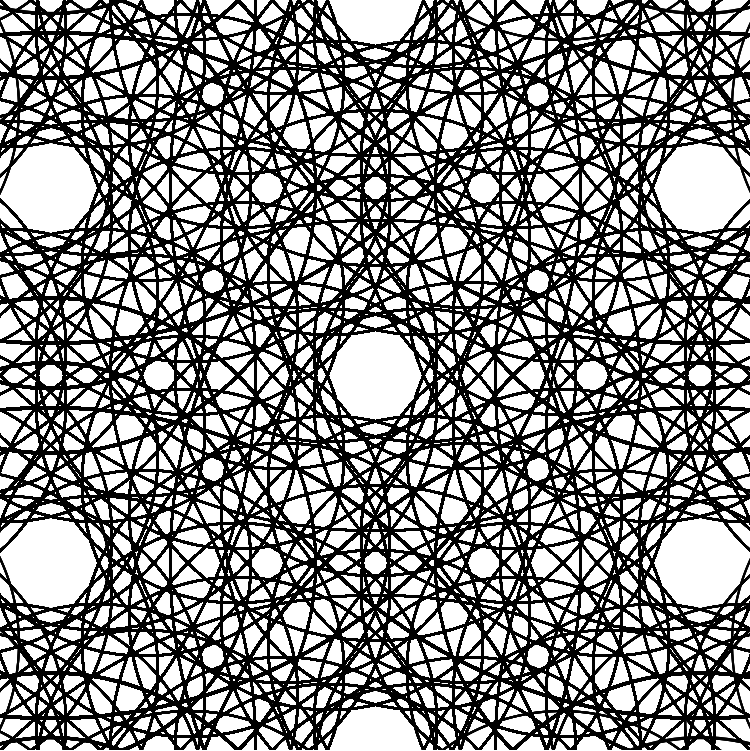}
\end{center}

\begin{center}
\includegraphics[width=10cm]{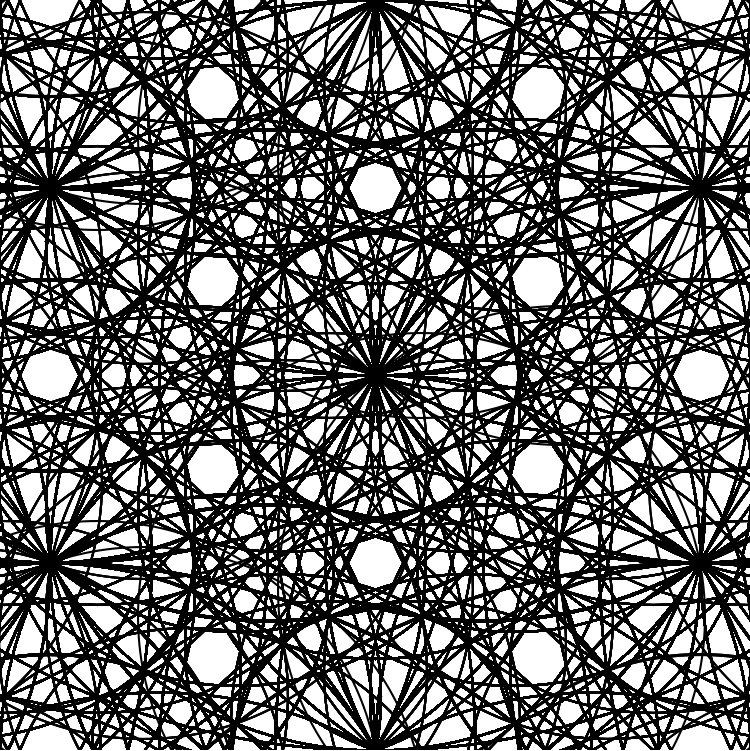}
\end{center}

The first part of Theorem \ref{theo:classification} in the
introduction concerns the classification up to conjugacy of the
maximal nonelementary $\CC$-Fuchsian subgroups of $\Ga_K$. Consider
the set $\maxcf$ of maximal nonelementary $\CC$-Fuchsian subgroups of
$\Ga_K$, on which the group $\Ga_K$ acts by conjugation. We will prove
that the discriminant map $\Ga\mapsto \Delta_\Ga$ on $\maxcf$ induces
a finite-to-one map from $\Ga_K\bs \maxcf$ onto $\NN-\{0\}$.

The second part of Theorem \ref{theo:classification} concerns the
classification up to commensurability and conjugacy. Given a group $G$
and a subgroup $H$ of $G$, recall that two subgroups $\Ga,\Ga'$ of $H$
are {\it commensurable} if $\Ga\cap\Ga'$ has finite index in $\Ga$ and
in $\Ga'$, and are {\it commensurable up to conjugacy in $G$} (or {\it
  commensurable in the wide sense}) if there exists $g\in G$ such that
$\Ga'$ and $g\Ga g^{-1}$ are commensurable.  For any positive natural
number $D$, let
$$
\Ga_{K,D}=\stab_{\Ga_K}[-D:0:1]\,.
$$
The group $\Ga_{K,D}$ is, by Proposition \ref{prop:stabmaxextCfuchs},
a maximal nonelementary extended $\CC$-Fuchsian subgroup, which
preserves the projective line $[-D:0:1]^\perp$. Its discriminant is
$2D$. We will prove that every element of $\maxcf$ with discriminant
$D$ is commensurable up to conjugacy in $\PSU_h$ with $\Ga_{K,2D}$.

\medskip
\noindent{\bf Proof of Theorem \ref{theo:classification}.} 
(1) Let $D\in\NN-\{0\}$ and let 
$$\maxcf(D)= \{\Ga\in
\maxcf\;:\;\Delta_\Ga=D\}\;.$$ Let 
$$
P_D=\left\{\begin{array}{l}{}[-\frac{D}{2}:0:1] {\rm ~if~} D 
{\rm ~is~even}\\
{}[0:1:0] {\rm ~if~} D=1\\
{}[-\frac{D-1}{2}:1:1] {\rm ~if~} D>1 {\rm ~is~odd}\;.
\end{array}\right.
$$
By Proposition \ref{prop:stabmaxextCfuchs}, the stabiliser in $\Ga_K$
of the positive rational point $P_D$ is a maximal nonelementary
extended $\CC$-Fuchsian subgroup of $\Ga_K$, with discriminant
$D$. Hence $\maxcf(D)$ is nonempty.

\medskip Let $\underline{G}$ be the semisimple connected linear
algebraic group defined over $\QQ$ such that $\underline{G}(\ZZ) =
\SU_h(\OOO_K)$ and $\underline{G}(\RR)=\SU_h$. Let $\pi: \underline{G}
\ra \GL(\underline{V})$ be the rational representation such that
$\underline{V}(\ZZ)=(\OOO_K)^3$, $\underline{V}(\RR)=\CC^3$ and
$\pi_{\mid \underline{G}(\RR)}$ is the linear action of $\SU_h$ on
$\CC^3$.  Let $\underline{X}_D$ be the closed algebraic submanifold of
$\underline{V}$ with equation $h=D$. In particular, $\underline{X}_D$
is defined over $\QQ$, and $\underline{X}_D(\RR)$ is homogeneous under
$\underline{G}(\RR) =\SU_h$, by Witt's theorem. The map 
$$
\underline{X}_D(\ZZ) = \underline{X}_D \cap \underline{V}(\ZZ)=
\{(z_0,z_1,z_2)\in (\OOO_K)^3\;:\;h(z_0,z_1,z_2)=D\}
\to \maxcf\,,
$$  which to $(z_0,z_1,z_2)$ associates the
stabiliser of $[z_0:z_1:z_2]$ in $\Ga_K$ (which is the image of
$\underline{G}(\ZZ)$ by the canonical map $\underline{G}(\RR)
=\SU_h\ra \PSU_h$), is well defined by Proposition
\ref{prop:stabmaxextCfuchs} and $\underline{G}(\ZZ)$-equivariant, and
its image contains $\maxcf(D)$. Hence the
finiteness of $\Ga_K\bs \maxcf(D)$ follows from
the finiteness of the number of orbits of $\underline{G}(\ZZ)$ on
$\underline{X}_D(\ZZ)$, see \cite[Thm.~6.9]{BorelHarishChandra62}.

\bigskip (2) Let $\Ga\in \maxcf$, and let
$D\in\NN-\{0\}$ be its discriminant. By Propositions
\ref{prop:Fuchsianrationalpoint} and \ref{prop:stabmaxextCfuchs}, and
by maximality, there is a unique positive rational point $P=
[z_0:z_1:z_2]$ with $z_0,z_1,z_2$ relatively prime in $\OOO_K$ such
that $\Ga$ is contained with finite index in $\stab_{\Ga_K}P$, and
$D=h(z_0,z_1,z_2)$.

\medskip\noindent{\bf Claim. } 
There exists $\ga\in\PSU_h(K)$ such that $\ga
P=[-2D:0:1]$. 

\medskip Assuming this claim for the moment, we conclude the proof of
the second part of Theorem \ref{theo:classification}: The groups
$\ga\Ga\ga^{-1}$ and $\Ga_{K,2D}$ are commensurable, since
$$
\ga\big(\stab_{\Ga_K}P\big)\ga^{-1}\cap \Ga_{K,2D}=\stab_{\ga\Ga_K\ga^{-1}\cap
  \Ga_K}\ga P =\ga\big(\stab_{\Ga_K\cap \ga^{-1}\Ga_K\ga}
P\big)\ga^{-1}
$$ 
and since $\PSU_h(K)$ is contained in the commensurator of
$\Ga_K=\PSU_h(\OOO_K)$ in $\PSU_h$ by a standard argument of
reduction to a common denominator.  \cqfd

\medskip The following result, useful for the proof of the above
claim, also gives a natural condition for when two such groups
$\Ga_{K,D}$ for $D\in\NN-\{0\}$ are commensurable up to conjugacy in
$\PSU_h$. A necessary and sufficient condition will be given as a
consequence of Section \ref{sect:quatalg}.

\blemm\label{lem:commensurablediscriminants} 
If $D,D'\in\NN-\{0\}$ satisfy $D'\in D \,N(\OOO_K)$, then $\Ga_{K,D}$
and $\Ga_{K,D'}$ are commensurable up to conjugacy in $\PSU_h(K)$.
\elemm

\dem Let $D\in\NN-\{0\}$ and $N\in N(\OOO_K)-\{0\}$. As seen above, we
only have to prove that there exists $\ga\in\PSU_h(K)$ such that
$\ga[-D:0:1]=[-DN:0:1]$.

Assume first that $D_K\equiv 0\mod 4$, so that $\OOO_K=\ZZ+
\frac{\sqrt{D_K}}{2} \ZZ$. Since $N\in N(\OOO_K)$, there exists
$x,y\in\ZZ$ such that $N=x^2-\frac{D_K}{4}y^2$. It is easy to check
using Equation \eqref{eq:calcstabpospoint} and since
$K=\QQ+i\sqrt{|D_K|}\QQ$ that the matrix
$$
\ga=\begin{pmatrix}  x & 0 & -\frac{i}{2}\sqrt{|D_K|}Dy\\ 0 & 1 & 0 \\
  -\frac{i}{2}\sqrt{|D_K|}\frac{y}{DN} & 0 & \frac{x}N \end{pmatrix}
$$ 
belongs to $\SU_h(K)$. Let $\ga$ be its image in $\PSU_h(K)$. It is
easy to check that as wanted $\ga[-D:0:1]=[-DN:0:1]$.

If $D_K\equiv 1\mod 4$, so that $\OOO_K=\ZZ+ \frac{1+\sqrt{D_K}}{2}
\ZZ$, the same argument works when $\ga$ in the above proof is
replaced by the matrix
$$
\begin{pmatrix}  x+\frac y2 & 0 & -\frac{i}{2}\sqrt{|D_K|}Dy\\ 0 & 1 & 0 \\
  -\frac{i}{2}\sqrt{|D_K|}\frac{y}{DN} & 0 & \frac{x+\frac y2}N \end{pmatrix}\,
$$
and the equation $N=x^2+xy+\frac{1-D_K}4 y^2$ with
$x,y\in\ZZ$.  \cqfd

\bigskip\noindent{\bf Proof of the claim. } As the lattice $\Ga_K$
does not preserve the complex geodesic with equation $z_2=0$, we may
assume that $z_2$ is nonzero, up to replacing $P$ by an element in its
orbit under $\Ga_K$, which does not change the discriminant $D$ of
$\Ga$.  Let $\ga_1$ be the Heisenberg translation by the element
$$
\Big[w_0=\frac{|z_1|^2}{2\,|z_2|^2} +i\,\Im\;\frac{z_0}{z_2}:
w=-\frac{z_1}{z_2}:1\Big]\in \Heis_3\;,
$$
which belongs to $\PSU_h(K)$. An easy computation shows that 
$$
\ga_1[z_0:z_1:z_2]=[-D:0:2\,N(z_2)]\;.
$$
Let $\ga_2$ be the image in $\PSU_h(K)$ of the diagonal element
$\begin{pmatrix} 2N(z_2) & 0 & 0\\ 0 & 1 & 0 \\ 0 & 0 &
  \frac{1}{2N(z_2)}\end{pmatrix}$ in $\SU_h(K)$. Then $\ga_2\ga_1$
maps $P$ to $[-2\,D\,N(z_2):0:1]$. By the previous lemma, there exists
$\ga_3\in\PSU_h(K)$ such that $\ga_3[-2\,D\,N(z_2):0:1]= [-2\,D:0:1]$.
Hence the claim follows with $\ga=\ga_3\ga_2 \ga_1$.  \cqfd

\section{Quaternion algebras}\label{sec:quaternion}
\label{sect:quatalg}

Let $a,b\in\ZZ$ with $a>0$ and $b<0$. The quaternion algebra $A=
\big(\frac{a,b}{\QQ}\big)$ is the $4$-dimensional central simple
algebra over $\QQ$ with standard generators $i,j,k$ satisfying the
relations $i^2=a$, $j^2=b$ and $ij=-ji=k$.  The {\em (reduced) norm}
of an element of $A$ is
$$
n(x_0+x_1i+x_2j+x_3k)=x_0^2-ax_1^2-bx_2^2+abx_3^2\;.
$$
The group of elements in $A(\ZZ)=\ZZ+i\ZZ+j\ZZ+k\ZZ$ with norm $1$ is
denoted by $A(\ZZ)^1$.  We refer to \cite{Vigneras80} and
\cite{MacRei03} for generalities on quaternion algebras.

\blemm
The map $\sigma=\sigma_{a,b}:A\to\mat_3(\CC)$ defined by 
$$
(x_0+x_1i+x_2j+x_3k)\mapsto\begin{pmatrix}
x_0+x_1\sqrt a & 0 & (x_2+x_3\sqrt a)\sqrt{b}\\ 0 & 1 & 0\\
(x_2-x_3\sqrt a)\sqrt{b} & 0 & x_0-x_1\sqrt a
\end{pmatrix}
$$
is a morphism of $\QQ$-algebras and $\sigma(A(\ZZ)^1)$ is a discrete
subgroup of the stabiliser of $[0:1:0]$ in $\SU_h$.  
\elemm

\dem It is well-known (and easy to check), see for instance
\cite{Katok92,MacRei03}, that the map $\sigma':A\to\mat_2(\RR)$
defined by
$$
(x_0+x_1i+x_2j+x_3k)\mapsto\begin{pmatrix}
x_0+x_1\sqrt a & (x_2+x_3\sqrt a)\sqrt{|b|}\\
-(x_2-x_3\sqrt a)\sqrt{|b|} & x_0-x_1\sqrt a
\end{pmatrix}\,
$$
is a morphism of $\QQ$-algebras and that the image of $A(\ZZ)^1$ is a
discrete subgroup of $\SLR$.  The map
$$
\iota:\begin{pmatrix}
a&b\\c&d
\end{pmatrix}\mapsto\begin{pmatrix}
a&0&ib\\0&1&0\\-ic&0& d
\end{pmatrix}
$$
is a morphism of $\QQ$-algebras, sending $\SLR$ into the stabiliser of
$[0:1:0]$ in $\SU_h$ (see Equation \eqref{eq:calcstabpospoint}). The
claim follows by noting that $\sigma=\iota\circ\sigma'$.  
\cqfd

\bigskip
\noindent{\bf Proof of Theorem \ref{theo:quaternion}.} By Theorem
\ref{theo:classification}, we only have to prove that the maximal
$\CC$-Fuchsian subgroup $F_D$ of $\Ga_K$ stabilising $[-2D:0:1]$
(which has finite index in $\Ga_{K,2D}$) arises from the quaternion
algebra $\big(\frac{D,\,D_K}{\QQ} \big)$. It is easy to check that the
element
$$
\ga_0=-\frac1{\sqrt 2}\begin{pmatrix} 
\sqrt D &\sqrt{2D}&\sqrt D\\
1 & 0 & -1\\
\frac1{2\sqrt{D}}& -\frac1{\sqrt{2D}}&\frac1{2\sqrt{D}}
\end{pmatrix}
$$
belongs to $\SU_h$ and maps $[0:1:0]$ to $[-2D:0:1]$.  Hence, using
Equation \eqref{eq:calcstabpospoint}, a matrix $M\in\SU_h(\OOO_K)$ has
its image (by the canonical projection $\SU_h\ra \PSU_h$) in $F_D$ if
and only if there exists $a,d\in\RR$ and $b,c\in i\RR$ with $ad-bc=1$
such that  $M=\ga_0 \begin{pmatrix} a & 0 & b \\ 0 & 1 & 0 \\ c & 0 & d 
\end{pmatrix}\ga_0^{-1}$.  A straightforward computation gives
$$
M=\begin{pmatrix} 
\frac{1}{4}(a+b+c+d+2) & \frac{\sqrt{D}}{2}(a-b+c-d) & 
\frac{D}{2}(a+b+c+d-2) \\ 
\frac{1}{4\sqrt{D}}(a+b-c-d) & \frac{1}{2}(a-b-c+d) & 
\frac{\sqrt{D}}{2}(a+b-c-d) \\
\frac{1}{8D}(a+b+c+d-2) & \frac{1}{4\sqrt{D}}(a-b+c-d) & 
\frac{1}{4}(a+b+c+d+2) 
\end{pmatrix}\;.
$$
This matrix has coefficients in $\OOO_K$ if and only if
$$
\begin{cases}
a+b+c+d-2\in 8D\OOO_K,\\
a+b-c-d\in 4\sqrt{D}\OOO_K,\\
a-b+c-d\in 4\sqrt{D}\OOO_K,\\
a-b-c+d\in 2\OOO_K\,.
\end{cases}
$$
Let $u=a+d$,
$v=\frac{1}{2\sqrt{D}}(a-d)$,  
$s'=b+c$ and $ t'=\frac{1}{2\sqrt{D}}(b-c)$. 
Hence $M$ has coefficients in $\OOO_K$ if
and only if
\begin{equation}\label{eq:equivuvst}
\begin{cases}
u+s'-2\in8D\OOO_K,\\ 
v+t'\in 2\OOO_K,\\
v-t'\in 2\OOO_K,\\
 u-s'\in 2\OOO_K\;.
\end{cases}
\end{equation}

Let $D'_K=\frac{D_K}{4}$ if $D_K\equiv 0 \mod 4$ and $D'_K=D_K$
otherwise. Recall that $\OOO_K\cap\RR=\ZZ$ and $\OOO_K\cap i\RR=
\ZZ\sqrt{D'_K}$.  
The equations \eqref{eq:equivuvst} imply in particular that
$u,v,s',t'\in\OOO_K$.  Note that $a,d\in\RR$ is equivalent to
$u,v\in\RR$, and $c,b\in i\RR$ is equivalent to $s',t'\in i\RR$.
Hence $u,v\in\ZZ$ and there exists $s,t\in\ZZ$ such that
$s'=s\sqrt{D'_K}, t'=t\sqrt{D'_K}$.
Therefore
$$
\ga_0^{-1}\,F_D\,\ga_0=\left\{
\begin{pmatrix} \frac u2+v \sqrt D&0&(\frac s2+t\sqrt D)\sqrt{D'_K}\\
0&1&0\\
(\frac s2-t\sqrt D)\sqrt{D'_K}&0& \frac u2-v \sqrt D
\end{pmatrix}:\begin{array}{l} u,v,s,t\in\ZZ\\
v+t\sqrt{D'_K}\in 2\OOO_K\\v-t\sqrt{D'_K}\in 2\OOO_K\\ 
u-s\sqrt{D'_K}\in2\OOO_K\\
u+s\sqrt{D'_K}-2\in 8D\OOO_K\end{array}\right\}\,.
$$

The group $\ga_0^{-1}\,F_D\,\ga_0$ is contained in
$\sigma_{D,D'_K}(A(\ZZ)^1)$, since the parameters $u$ and $s$ have to
be even as a consequence of the defining equations of
$\ga_0^{-1}\,F_D\,\ga_0$. Furthermore, $\ga_0^{-1}\,F_D\,\ga_0$
contains $\sigma_{D,D'_K}(\OOO^1)$, where $\OOO$ is the order of $A$
defined by
$$
\OOO=\{x_0+ix_1+jx_2+kx_3\in A(\ZZ)\;:\; 
x_1,x_2,x_3\equiv 0 \mod 4D\}\;.
$$
Indeed, if $x_0+ix_1+jx_2+kx_3\in\OOO^1$, then with $u=2x_0$,
$s=2x_2$, $v=x_3$, $t=x_4$, we have, since $x_0\equiv 1\mod 4D$ by the
condition $n(x_0+ix_1+jx_2+kx_3)=1$,
$$
\begin{cases}
v\pm t\sqrt{D'_K}\in 2\ZZ +2\sqrt{D'_K}\,\ZZ\subset 2\OOO_K\\
u-s\sqrt{D'_K}\in 2\ZZ +2\sqrt{D'_K}\,\ZZ\subset 2\OOO_K\\
u-2+s\sqrt{D'_K}=2(x_0-1)+2x_2\sqrt{D'_K}\in 
8D\ZZ +8D\sqrt{D'_K}\,\ZZ\subset 8D\OOO_K\;.
\end{cases}
$$
Since $\sigma_{D,D'_K}(\OOO^1)$ has finite index in
$\sigma_{D,D'_K}(A(\ZZ)^1)$ (see for instance
\cite{Vigneras80}, Coro.~1.5 in Chapt.~IV), the groups $\ga_0^{-1}\,F_D\,\ga_0$ and
$\sigma_{D,D'_K}(A(\ZZ)^1)$ are hence commensurable.

Since $\big(\frac{D,D'_K}{\QQ}\big)=\big(\frac{D,D_K}{\QQ}\big)$ as
$D'_K$ and $D_K$ differ by a square, the result follows. 
\cqfd

\bigskip
In particular, a maximal nonelementary $\CC$-Fuchsian subgroup of
$\Ga_K$ of discriminant $D$ is cocompact (in its copy of
$\PSL_2(\RR)$) if and only if $\big(\frac{D,\,D_K}{\QQ} \big)$ is a
division algebra (see for instance \cite[Thm.~5.4.1]{Katok92}).

\bigskip 
The following corollaries follow from the arguments in
\cite{Maclachlan86}, pages 309 and 310.  Corollary \ref{coro:intro} of
the introduction follows from Corollary \ref{coro:manyclasses} below.

\bprop
Let $A$ be an indefinite quaternion algebra over $\QQ$. There
exists an arithmetic $\CC$-Fuchsian subgroup of $\Ga_K$ whose
associated quaternion algebra is $A$ if and only if the primes  at
which $A$ is ramified are either ramified or inert in $K$. \cqfd   
\eprop

\bcoro\label{coro:manyclasses} (Chinburg-Stover) Every Picard modular group
$\SU_h(\OOO_K)$ contains infinitely many wide commensurability classes
of maximal nonelementary $\CC$-Fuchsian subgroups. \cqfd 
\ecoro

\bcoro 
Any arithmetic Fuchsian group whose associated quaternion algebra is
defined over $\QQ$ is contained (up to commensurability) as a
$\CC$-Fuchsian subgroup of some Picard modular group $\Ga_K$. \cqfd
\ecoro

\bcoro For all quadratic irrational number fields $K$ and $K'$, there
are infinitely many commensurability classes of arithmetic Fuchsian
subgroups with representatives in both Picard modular groups $\Ga_K$
and $\Ga_{K'}$. \cqfd 
\ecoro

{\small \bibliography{../biblio} }

\bigskip
{\small
\noindent \begin{tabular}{l} 
Department of Mathematics and Statistics, P.O. Box 35\\ 
40014 University of Jyv\"askyl\"a, FINLAND.\\
{\it e-mail: jouni.t.parkkonen@jyu.fi}
\end{tabular}
\medskip

\noindent \begin{tabular}{l}
D\'epartement de math\'ematique, UMR 8628 CNRS, B\^at.~425\\
Universit\'e Paris-Sud,
91405 ORSAY Cedex, FRANCE\\
{\it e-mail: frederic.paulin@math.u-psud.fr}
\end{tabular}
}

\end{document}